\documentclass[twoside,twocolumn,9pt]{article}
\usepackage{extsizes}
\usepackage[super,sort&compress,comma]{natbib} 
\usepackage[version=3]{mhchem}
\usepackage[left=1.5cm, right=1.5cm, top=1.785cm, bottom=2.0cm]{geometry}
\usepackage{balance}
\usepackage{times,mathptmx}
\usepackage{sectsty}
\usepackage{graphicx}
\usepackage{lastpage}
\usepackage[format=plain,justification=raggedright,singlelinecheck=false,font={stretch=1.125,small,sf},labelfont=bf,labelsep=space]{caption}
\usepackage{float}
\usepackage{fancyhdr}
\usepackage{fnpos}
\usepackage[english]{babel}
\usepackage{array}
\usepackage{droidsans}
\usepackage{charter}
\usepackage[T1]{fontenc}
\usepackage[usenames,dvipsnames]{xcolor}
\usepackage{setspace}
\usepackage[compact]{titlesec}

\usepackage{epstopdf}

\definecolor{cream}{RGB}{222,217,201}

\usepackage{amsmath,amssymb} 
\renewcommand{\P}{\mathbb{P}}
\newcommand{\E}{\mathbb{E}}
\newcommand{\N}{\mathbb{N}}
\newcommand{\R}{\mathbb{R}}
\renewcommand{\div}{{\rm div}}

\begin{document}

%

\makeFNbottom
\makeatletter
\renewcommand\LARGE{\@setfontsize\LARGE{15pt}{17}}
\renewcommand\Large{\@setfontsize\Large{12pt}{14}}
\renewcommand\large{\@setfontsize\large{10pt}{12}}
\renewcommand\footnotesize{\@setfontsize\footnotesize{7pt}{10}}
\makeatother

\renewcommand{\thefootnote}{\fnsymbol{footnote}}
\renewcommand\footnoterule{\vspace*{1pt}%
\color{cream}\hrule width 3.5in height 0.4pt \color{black}\vspace*{5pt}} 
\setcounter{secnumdepth}{5}

\makeatletter 
\renewcommand\@biblabel[1]{#1}            
\renewcommand\@makefntext[1]%
{\noindent\makebox[0pt][r]{\@thefnmark\,}#1}
\makeatother 
\renewcommand{\figurename}{\small{Fig.}~}
\sectionfont{\sffamily\Large}
\subsectionfont{\normalsize}
\subsubsectionfont{\bf}
\setstretch{1.125} 
\setlength{\skip\footins}{0.8cm}
\setlength{\footnotesep}{0.25cm}
\setlength{\jot}{10pt}
\titlespacing*{\section}{0pt}{4pt}{4pt}
\titlespacing*{\subsection}{0pt}{15pt}{1pt}

\fancyfoot{}
\fancyfoot[LO,RE]{\vspace{-7.1pt}\includegraphics[height=9pt]{head_foot/LF}}
\fancyfoot[CO]{\vspace{-7.1pt}\hspace{13.2cm}\includegraphics{head_foot/RF}}
\fancyfoot[CE]{\vspace{-7.2pt}\hspace{-14.2cm}\includegraphics{head_foot/RF}}
\fancyfoot[RO]{\footnotesize{\sffamily{1--\pageref{LastPage} ~\textbar  \hspace{2pt}\thepage}}}
\fancyfoot[LE]{\footnotesize{\sffamily{\thepage~\textbar\hspace{3.45cm} 1--\pageref{LastPage}}}}
\fancyhead{}
\renewcommand{\headrulewidth}{0pt} 
\renewcommand{\footrulewidth}{0pt}
\setlength{\arrayrulewidth}{1pt}
\setlength{\columnsep}{6.5mm}
\setlength\bibsep{1pt}

\makeatletter 
\newlength{\figrulesep} 
\setlength{\figrulesep}{0.5\textfloatsep} 

\newcommand{\topfigrule}{\vspace*{-1pt}%
\noindent{\color{cream}\rule[-\figrulesep]{\columnwidth}{1.5pt}} }

\newcommand{\botfigrule}{\vspace*{-2pt}%
\noindent{\color{cream}\rule[\figrulesep]{\columnwidth}{1.5pt}} }

\newcommand{\dblfigrule}{\vspace*{-1pt}%
\noindent{\color{cream}\rule[-\figrulesep]{\textwidth}{1.5pt}} }

\makeatother

\twocolumn[
  \begin{@twocolumnfalse}
\vspace{3cm}
\sffamily
\begin{tabular}{m{1.5cm} p{13.5cm} }

 &    \noindent\LARGE{\textbf{Jump Markov
    models and transition state theory: the Quasi-Stationary Distribution approach}} \\
\vspace{0.3cm} & \vspace{0.3cm} \\

 & \noindent\large{Giacomo Di Ges\`u\textit{$^{a}$}, Tony
   Leli\`evre\textit{$^{\ast a}$}, Dorian Le Peutrec\textit{$^{a,b}$} and Boris Nectoux\textit{$^{a}$}} \\

\

 &   \noindent\normalsize{We are
  interested in the connection between a metastable
  continuous state space Markov process (satisfying {\em e.g.} the Langevin or overdamped Langevin equation) and a jump Markov
  process in a discrete state space. More precisely, we use the notion of quasi-stationary
  distribution within a metastable state for the continuous state space
  Markov process to parametrize the exit event
  from the  state. This approach is useful to analyze and justify methods which use the jump
  Markov process underlying a metastable dynamics as a support to
  efficiently sample the state-to-state dynamics (accelerated dynamics
  techniques). Moreover, it is possible by this approach to quantify the error on the
  exit event when the parametrization of the jump Markov model is
  based on the  Eyring-Kramers formula. This therefore provides a
  mathematical framework to justify the use of transition state theory
  and the Eyring-Kramers formula to build kinetic Monte Carlo or Markov
  state models.}
\\

\

\end{tabular}

 \end{@twocolumnfalse} \vspace{0.6cm}

 ] 

\renewcommand*\rmdefault{bch}\normalfont\upshape
\rmfamily
\section*{}
\vspace{-1cm}


\footnotetext{\textit{$^{a}$~CERMICS, \'Ecole des Ponts, Universit\'e
    Paris-Est, INRIA, 77455 Champs-sur-Marne, France. E-mail: \{di-gesug,lelievre,nectoux\}@cermics.enpc.fr}}
\footnotetext{\textit{$^{b}$~Laboratoire de Math\'ematiques d'Orsay,
    Univ. Paris-Sud, CNRS, Université Paris-Saclay, 91405 Orsay, France. E-mail: dorian.lepeutrec@math.u-psud.fr}}


\footnotetext{$^\ast$~Corresponding author. This work is supported by the European Research Council under the
European Union's Seventh Framework Programme (FP/2007-2013) / ERC
Grant Agreement number 614492.}





\section{Introduction and motivation}

Many theoretical studies and numerical methods in materials science~\cite{voter-05},
biology~\cite{bowman-pande-noe-14} and chemistry, aim at modelling the dynamics at the atomic
level as a jump Markov process between states. Our objective in this
article is to discuss the relationship between such a mesoscopic model
(a Markov process over a discrete state space)
and the standard microscopic full-atom description (typically a Markov
process over a continuous state space, namely a molecular dynamics simulation).

The objectives of a modelling using a jump Markov process rather than a
detailed microscopic description at the atomic level are
numerous. From a modelling viewpoint, new insights can be
gained by building coarse-grained models, that are easier to handle. From a numerical viewpoint,
the hope is to be able to build the jump Markov process from short
simulations of the full-atom dynamics. Moreover, once the
parametrization is done, it is possible to simulate the system over
much longer timescales than the time horizons attained by standard
molecular dynamics, either by using directly the jump Markov process,
or as a support to accelerate molecular
dynamics~\cite{voter-97,voter-98,sorensen-voter-00}. It is also
possible to use dedicated algorithms to extract from the graph
associated with the jump Markov process the most important features of
the dynamics (for example quasi-invariant sets and essential
timescales using large deviation theory~\cite{freidlin-wentzell-84}), see for example~\cite{wales-03,cameron-2014}.

In order to parametrize the jump Markov process, one needs to define
rates from one state to another. The concept of jump rate between two
states is one of
the fundamental notions in the modelling of materials. Many papers
have been devoted to the rigorous evaluation of jump rates from a full-atom description. The most
famous formula is probably the rate derived in the harmonic transition
state
theory~\cite{marcelin-15,polanyi-eyring-31,eyring-35,wigner-38,horiuti-38,kramers-40,vineyard-57},
which gives an explicit expression for the rate in terms of the
underlying potential energy function (see the Eyring-Kramers formula~\eqref{eq:EK}
below). See for example the review
paper\cite{hanggi-talkner-barkovec-90}.

Let us now present the two models: the jump Markov model, and the
full-atom model, before discussing how the latter can be related to
the former.

\subsection{Jump Markov models}\label{sec:Z}

Jump Markov models are continuous-time Markov processes with values
in a discrete state space. In the context of molecular modelling,
they are known as Markov state models~\cite{bowman-pande-noe-14,schuette-sarich-13} or
kinetic Monte Carlo models~\cite{voter-05}. They consist of a
collection of states that we can assume to be indexed by integers, and rates $(k_{i,j})_{i \neq j
  \in \N}$ which are associated with transitions between these
states. For a state $i \in \N$, the states $j$ such
that $k_{i,j} \neq 0$ are the neighboring states of $i$ denoted in the
following by
\begin{equation}\label{eq:Ni}
{\mathcal N}_i=\{j \in \N, \, k_{i,j} \neq 0\}.
\end{equation}
 One can thus think of a jump Markov model as a graph: the states are the vertices,
and an oriented edge between two vertices $i$ and $j$ indicates that
$k_{i,j}\neq 0$.

Starting at time $0$ from a state $Y_0 \in \N$, the model consists in iterating the following two steps
over $n \in \N$: Given $Y_n$,
\begin{itemize}
\item Sample the residence time $T_n$ in $Y_n$ as an exponential random
  variable with parameter $\sum_{j \in {\mathcal N}_{Y_n}} k_{Y_n,j}$:
\begin{equation}\label{eq:T_n}
\forall t \ge 0, \, \P(T_n\ge t | Y_n=i) = \exp\left(- \left[\sum_{j
      \in {\mathcal N}_{i}} k_{i,j}\right] \, t\right).
\end{equation}
\item Sample independently from $T_n$ the next visited state $Y_{n+1}$ starting from $Y_n$
  using the following law
\begin{equation}\label{eq:Y_n}
\forall j  \in {\mathcal N}_{i}, \, \P(Y_{n+1}=j | Y_n=i) = \frac{k_{i,j}}{\sum_{j
    \in {\mathcal N}_{i}} k_{i,j}}.
\end{equation}
\end{itemize}
The associated continuous-time process $(Z_t)_{t \ge 0}$
with values in $\N$ defined by:
\begin{equation}\label{eq:Z}
\forall n \ge 0, \, \forall t \in \left[\sum_{m=0}^{n-1} T_m,
\sum_{m=0}^{n} T_m\right), \quad Z_t = Y_n
\end{equation}
(with the convention $\sum_{m=0}^{-1}=0$) is then a (continous-time) jump Markov process.

\subsection{Microscopic dynamics}\label{sec:LOL}

At the atomic level, the basic ingredient is a potential energy
function $V:\R^d \to \R$ which to a set of positions of atoms in
$x \in \R^d$ (the dimension $d$ is typically 3 times the number of atoms)
associates an energy $V(x)$. In all the following, we assume that $V$
is a smooth Morse function: for each $x \in \R^d$, if $x$ is a
critical point of $V$ (namely if $\nabla
V(x)=0$), then the Hessian $\nabla^2 V(x)$ of $V$ at point $x$ is a
nonsingular matrix. From this function $V$, dynamics are
built such as the Langevin dynamics:
\begin{equation}\label{eq:L}
\begin{aligned}
dq_t &= M^{-1} p_t \, dt\\
dp_t &= -\nabla V(q_t) \, dt - \gamma M^{-1} p_t \, dt + \sqrt{2
  \gamma \beta^{-1}} dW_t
\end{aligned}
\end{equation}
or the overdamped Langevin dynamics:
\begin{equation}\label{eq:OL}
dX_t = -\nabla V(X_t) \, dt + \sqrt{2 \beta^{-1}} dW_t.
\end{equation}
Here, $M \in \R^{d\times d}$ is the mass matrix, $\gamma >0$ is the
friction parameter, $\beta^{-1}=k_B T > 0$ is the inverse
temperature and $W_t \in \R^d$ is a $d$-dimensional Brownian motion. The Langevin dynamics gives the evolution of the positions
$q_t \in \R^d$ and the momenta $p_t \in \R^d$. The overdamped Langevin
dynamics is in position space: $X_t \in \R^d$. The overdamped Langevin
dynamics is derived from the Langevin dynamics in the large friction
limit and using a rescaling in time: assuming $M={\rm Id}$ for simplicity,
in the limit $\gamma
\to \infty$, $(q_{\gamma t})_{t \ge 0}$
converges to $(X_t)_{t \ge 0}$ (see for example Section 2.2.4 in~\cite{lelievre-rousset-stoltz-book-10}).

\subsection{From a microscopic dynamics to a jump Markov dynamics}

Let us now discuss how one can relate the  microscopic dynamics~\eqref{eq:L}
or~\eqref{eq:OL} to the jump Markov model~\eqref{eq:Z}. The basic
observation which justifies why this question is relevant is the
following. It is observed that, for applications in biology, material
sciences or chemistry, the microscopic dynamics~\eqref{eq:L}
or~\eqref{eq:OL} are metastable. This means that the stochastic
processes $(q_t)_{t \ge 0}$ or $(X_t)_{t \ge 0}$ remain trapped for a
long time in some region of the configurational space (called a
metastable region) before hopping to another metastable
region. Because the system remains for very long times in a metastable
region before exiting, the hope is that it loses the memory of the way
it enters, so that the exit event from this
region can be modelled as one move of a jump Markov process such as~\eqref{eq:Z}.

\begin{figure}[h]
\centering
  \includegraphics[height=5cm]{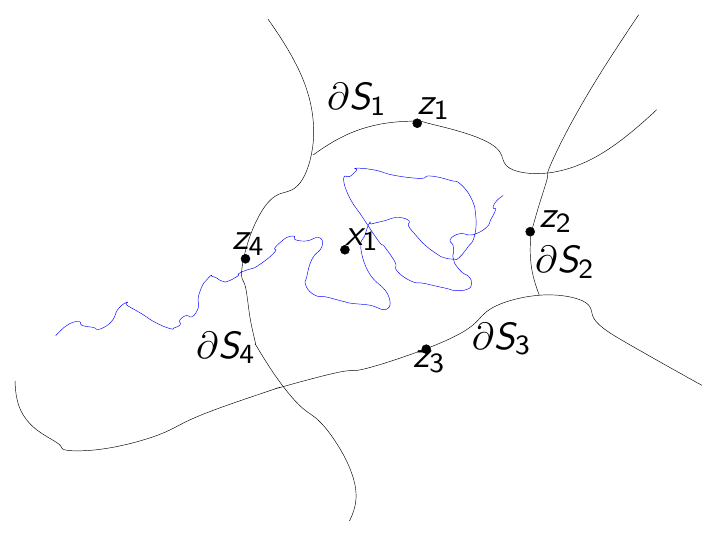}
  \caption{The domain $S$. The boundary $\partial S$ is divided into
    $4$ subdomains $(\partial S_i)_{1 \le i \le 4}$, which are the
    common boundaries with the neighboring states.}
  \label{fig:S}
\end{figure}

Let us now consider a subset $S \subset \R^d$ of the configurational
space for the microscopic dynamics. Positions in $S$ are associated with one of the discrete state in
$\N$ of $(Z_t)_{t \ge 0}$, say the state
$0$ without loss of generality. If $S$ is
metastable (in a sense to be made precise), it should be possible to
justify the fact that the exit event can be modeled using a jump
Markov process, and to compute the associated exit rates $(k_{0,j})_{j
  \in {\mathcal N}_0}$
from the state $0$ to the neighboring states using the
dynamics~\eqref{eq:L} or~\eqref{eq:OL}. The aim of this paper
is precisely to discuss these questions and in particular to
prove rigorously under which assumption the Eyring-Kramers
formula can be used to estimate the exit rates $(k_{0,j})_{j \in
  {\mathcal N}_0}$,
namely: 
\begin{equation}\label{eq:EK}
\forall j \in {\mathcal N}_0, \, k_{0,j}=\nu_{0,j} \exp(-\beta[V(z_j)-V(x_1)])
\end{equation}
where $\nu_{0,j}>0$ is a prefactor, $x_1=\arg\min_{x \in S} V(x)$ and
$z_j=\arg\min_{z \in \partial S_j} V(z)$ where $\partial S_j
\subset \partial S$ denotes the part of the boundary $\partial S$
which connects the state $S$ (numbered $0$) with the subset of $\R^d$
associated with state numbered $j \in {\mathcal N}_0$. See Figure~\ref{fig:S}.

The prefactor $\nu_{0,j}$ depends on the dynamic under
consideration and on $V$ around $x_1$ and $z_j$. Let us give a few
examples. If $S$ is taken as the basin of attraction of $x_1$ for the
dynamics $\dot{x}=-\nabla V(x)$ so that the points $z_j$ are order
one saddle points, the prefactor writes for the Langevin
dynamics~\eqref{eq:L} (assuming again $M={\rm Id}$ for simplicity):
\begin{equation}\label{eq:nu_L}
\nu^{L}_{0,j}=\frac{1}{4\pi} \left(\sqrt{\gamma^2 +4|\lambda^-(z_j)|}
  - \gamma \right)
\frac{\displaystyle \sqrt{\det (\nabla^2 V)(x_1)}}{\displaystyle \sqrt{|\det (\nabla^2 V)(z_j)|}}
\end{equation}
where, we recall, $\nabla^2 V$ is the Hessian of $V$, and $\lambda^-(z_j)$ denotes
the negative eigenvalue of $\nabla^2 V(z_j)$. This
formula was derived by Kramers in~\cite{kramers-40} in a
one-dimensional situation.
The equivalent formula for the overdamped Langevin
dynamics~\eqref{eq:OL} is:
\begin{equation}\label{eq:nu_OL}
\nu^{OL}_{0,j}=\frac{1}{2\pi}|\lambda^-(z_j)| \frac{\displaystyle \sqrt{\det (\nabla^2 V)(x_1)}}{\displaystyle\sqrt{|\det (\nabla^2 V)(z_j)|}}.
\end{equation}
Notice that $\lim_{\gamma \to \infty} \gamma \nu^L_{0,j}=
\nu^{OL}_{0,j}$, as expected from the rescaling in time used to go
from Langevin to overdamped Langevin (see Section~\ref{sec:LOL}). The
formula~\eqref{eq:nu_OL} has again been obtained by Kramers
in~\cite{kramers-40}, but also by many authors previously,
see the exhaustive review of the literature reported
in~\cite{hanggi-talkner-barkovec-90}. In Section~\ref{sec:EK_litt} below, we will
review mathematical results where formula
\eqref{eq:nu_L}--\eqref{eq:nu_OL} are rigorously derived.

In practice, there are thus two possible approaches to determine the
rates $(k_{i,j})$. On the one hand, when the number of states is not too large, one can
precisely study the transitions between metastable states for the
microscopic dynamics using dedicated algorithms~\cite{henkelman-johanneson-jonsson-02,e-vanden-eijnden-10}: the nudged elastic
band~\cite{jonsson-mills-jacobsen-98}, the string
method~\cite{e-ren-vanden-eijnden-02,e-ren-vanden-eijnden-05} and the
max flux approach~\cite{zhao-shen-skeel-10} aim at finding one typical
representative path. Transition path sampling
methods~\cite{dellago-bolhuis-chandler-99,dellago-bolhuis-09} sample
transition paths starting from an initial guessed trajectory, and using a Metropolis
Hastings algorithm in path space. Other approaches aim at sampling the
ensemble of paths joining two metastable states, without any initial
guess: see the Adaptive Multilevel Splitting method~\cite{cerou-guyader-07a,cerou-guyader-lelievre-pommier-11},
transition interface
sampling~\cite{van-erp-moroni-bolhuis-03,van-erp-bolhuis-05}, forward
flux
sampling~\cite{allen-warren-ten-wolde-05,allen-valeriani-ten-wolde-09},
milestoning
techniques~\cite{faradjian-elber-04,maragliano-vanden-eijnden-roux-09,schuette-noe-lu-sarich-vanden-einjden-11}
and the associated Transition Path Theory which gives a nice
mathematical
framework~\cite{e-vanden-eijnden-04,vanden-eijnden-venturoli-ciccotti-elber-08,lu-nolen-15}. On
the other hand, if
the number of states is very large, it may be too cumbersome to sample
all the transition paths, and one can use instead the Eyring-Kramers formula~\eqref{eq:EK}, which
requires to look for the local minima and the order one saddle points of
$V$, see for example~\cite{wales-03}. Algorithms
to look for saddle points include the dimer method
\cite{henkelman-jonsson-99,zhang-du-12}, activation relaxation
techniques \cite{barkema-mousseau-96,mousseau-barkema-98}, or the
gentlest ascent dynamics \cite{samanta-e-12}, for example.

The aim of this paper is threefold. First, we would like to give a
mathematical setting to quantify the metastability of a domain $S
\subset \R^d$ for a microscopic dynamics such as~\eqref{eq:L} or~\eqref{eq:OL},
and to rigorously justify the fact that for a metastable domain, the
exit event can be modeled using a jump process such as~\eqref{eq:Z}. This question is
addressed in Section~\ref{sec:QSD}, where we introduce the notion of
quasi-stationary distribution. Second, we explain in
Section~\ref{sec:Voter} how this framework can be used to analyze
algorithms which have been proposed by A.F. Voter, to accelerate the
sampling of the state-to-state dynamics using the underlying jump
Markov process. We will discuss in particular
the Parallel Replica algorithm~\cite{voter-98}. Both these aspects
were already presented by the second author in previous works, see for
example the review paper~\cite{lelievre-15}. The main novelty of this article is in
Section~\ref{sec:EK}, which is devoted to a justification of the use of the
Eyring-Kramers formula~\eqref{eq:EK} in order to
parametrize jump Markov models.

Before getting to the heart of the matter, let us make two
preliminary remarks. First, the objective of this paper is to give a
self-contained overview of the interest of using the quasi-stationary
distribution to analyze metastable processes. For the sake of
conciseness, we therefore do not
provide extensive proofs of the results we present, but we 
give the relevant references when necessary. Second, we will concentrate in
the following on the overdamped Langevin dynamics~\eqref{eq:OL} when
presenting mathematical results. All the algorithms presented below
equally apply to (and are actually used on) the Langevin
dynamics~\eqref{eq:L}. As will be explained below, the notion of quasi-stationary distribution which is the cornerstone of our analysis is
also well defined for Langevin dynamics. However, the mathematical
analysis of Section~\ref{sec:EK} is for the moment restricted to the overdamped Langevin dynamics~\eqref{eq:OL}.


\section{Metastable state and quasi-stationary distribution}\label{sec:QSD}

The setting in this section is the following. We consider the
overdamped Langevin dynamics~\eqref{eq:OL} for simplicity\footnote{The existence of the QSD and the convergence of the
conditioned process towards the QSD for the Langevin
process \eqref{eq:L} follows from the recent paper~\cite{nier-14}.}
 and a subset
$S \subset \R^d$ which is assumed to be bounded and smooth. We would like to first characterize the fact that
$S$ is a metastable region for the
dynamics. Roughly speaking, metastability means that the local equilibration
time within $S$ is much smaller than the exit time from $S$. In order
to approximate the original dynamics by a jump Markov model, we need
such a separation of timescales (see the discussion in
Section~\ref{sec:ParRep_conc} on how to take into account non-Markovian features). Our first
task is to give a precise meaning to that. Then, if $S$ is metastable,
we would like to study the exit event from $S$, namely the exit time
and the exit point from $S$, and to see if it can be related to the
exit event for a jump Markov model (see~\eqref{eq:T_n}--\eqref{eq:Y_n}). The analysis will use the notion of quasi-stationary distribution (QSD), that we now introduce.

\subsection{Definition of the QSD}\label{sec:QSD_def}
Consider the first exit
time from $S$:
\[\tau_S=\inf \{t \ge 0, \, X_t \not\in S \},\]
where $(X_t)_{t \ge 0}$ follows the
overdamped Langevin dynamics~\eqref{eq:OL}.

A probability measure $\nu_S$ with support in $S$ is called a QSD for
the Markov process $(X_t)_{t \ge 0}$ if and only if 
\begin{equation}\label{eq:QSD}
\nu_S(A)=\frac{{\displaystyle\int_S \P^x(X_t
 \in A, t < \tau_S) \, \nu_S(d x)}}{{\displaystyle\int_S \P^x(t < \tau_S) \,
  \nu_S(d x)}}, \quad 
\forall t > 0, \, \forall A \subset S.
\end{equation}

Here and in the following, $\P^x$ denotes the probability measure under which $X_0=x$.
In other words, $\nu_S$ is a QSD if, when $X_0$ is distributed according
to $\nu_S$, the law of $X_t$, 
conditional on $(X_s)_{0
  \le s \le t}$ remaining in the state $S$, is still $\nu_S$, for all
positive~$t$.

The QSD satisfies three properties which will be crucial in the
following. We refer for example to \cite{le-bris-lelievre-luskin-perez-12} for detailed proofs of these results and to \cite{collet-martinez-san-martin-13} for more
general results on QSDs.

\subsection{First property: definition of a metastable state}\label{sec:QSD1}
Let $(X_t)_{t \ge 0}$ follow the dynamics \eqref{eq:OL}
with an initial condition $X_0$ distributed according to a
distribution $\mu_0$ with support in $S$. Then there exists a probability
distribution $\nu_S$ with support in $S$ such that, for any initial
distribution $\mu_0$  with support in $S$,
\begin{equation}\label{eq:QSD_1}
\lim_{t \to \infty} {\rm Law}(X_t | \tau_S > t) = \nu_S.
\end{equation}
The distribution $\nu_S$ is the QSD associated with $S$.

A consequence of this proposition is the existence and uniqueness of
the QSD. The QSD is the long-time limit of the law of the (time
marginal of the) process
conditioned to stay in the state~$S$: it can be seen as a `local ergodic
measure' for the stochastic process in $S$. 

This proposition gives a first intuition to properly define
a metastable state. A metastable state is a state such that the typical exit time is
much larger than the local equilibration time, namely the time to
observe the convergence to the QSD in \eqref{eq:QSD_1}. We will explain
below how to quantify this timescale discrepancy (see~\eqref{eq:S_metastable}) by identifying the rate of convergence
in~\eqref{eq:QSD_1} (see~\eqref{eq:erreur}).

\subsection{Second property: eigenvalue problem}\label{sec:QSD2}
Let $L=-\nabla V \cdot \nabla + \beta^{-1} \Delta$ be the
infinitesimal generator of $(X_t)_{t \ge 0}$ (satisfying \eqref{eq:OL}). Let us consider the
first eigenvalue and eigenfunction associated with the adjoint operator
$L^\dagger=\div ( \nabla V + \beta^{-1} \nabla)$ with 
homogeneous Dirichlet
boundary condition on $\partial S$:
\begin{equation}\label{eq:QSD_2}
\left\{
\begin{aligned}
L^\dagger u_1 &= - \lambda_1 u_1 &&\text{on $S$},\\
u_1 &=0 &&\text{on $\partial S$}.
\end{aligned}
\right.
\end{equation}
Then, the QSD $\nu_S$ associated with $S$ satisfies
\[d \nu_S = \frac{u_1(x)  \, d x}{\displaystyle\int_S u_1(x)  \, d x}  \]
where $d x$ denotes the Lebesgue measure on $S$.

Notice that $L^\dagger$ is a negative operator in $ L^2(e^{\beta V}) $  so that $\lambda_1 >
0$. Moreover, it follows from general results on the first eigenfunction of elliptic
operators that $u_1$ has a sign on $S$, so that one can choose without
loss of generality $u_1>0$.

The QSD thus has a density with respect to 
Lebesgue measure, which
is simply the ground state of the Fokker--Planck operator $L^\dagger$
associated with the dynamics with absorbing boundary conditions. This
will be crucial in order to analyze the Eyring-Kramers formula in Section~\ref{sec:EK}.

\subsection{Third property: the exit event}\label{sec:QSD3}
Finally, the third property of the QSD concerns the exit event
starting from the QSD.
Let us assume that $X_0$ is distributed according to the QSD $\nu_S$ in $S$. Then the
law of the pair
 $(\tau_S,X_{\tau_S})$ (the first exit time and the
first exit point) is fully characterized by the following properties:
(i) $\tau_S$ is exponentially distributed with parameter $\lambda_1$
(defined in \eqref{eq:QSD_2}); (ii) $\tau_S$ is independent of
$X_{\tau_S}$; (iii) The law of $X_{\tau_S}$ is the following: for any bounded measurable
  function $\varphi: \partial S \to \R$,
\begin{equation}\label{eq:QSD_3}
\E^{\nu_S}(\varphi(X_{\tau_S})) = - \frac{\displaystyle\int_{\partial S}
  \varphi \, \partial_n u_1  \, d  \sigma}{\displaystyle\beta \lambda_1 \int_{S}
  u_1(x)  \, d x} ,
\end{equation}
where $\sigma$ denotes the Lebesgue measure on $\partial S$ and
$\partial_n u_1=\nabla u_1 \cdot n$ denotes the outward normal
derivative of $u_1$ (defined in \eqref{eq:QSD_2}) on $\partial S$. The
superscript $\nu_S$ in $\E^{\nu_S}$ indicates that the
initial condition $X_0$ is assumed to be distributed according to
$\nu_S$. 

\subsection{Error estimate on the exit event}\label{sec:erreur}

We can now state a result concerning the error made when approximating
the exit event of the process which remains for a long time in $S$ by
the exit event of the process starting from the QSD. The following
result is proven in~\cite{le-bris-lelievre-luskin-perez-12}.
Let $(X_t)_{t \ge 0}$ satisfy \eqref{eq:OL} with $X_0
\in S$. Introduce the first  two eigenvalues
$-\lambda_2 < -\lambda_1 < 0$
of the operator $L^\dagger$ on $S$ with homogeneous Dirichlet boundary
conditions on $\partial S$ (see Section~\ref{sec:QSD2}). Then there
exists a constant $C>0$ (which depends on the law of $X_0$), such that,
for all $t \ge \frac{C}{(\lambda_2  - \lambda_1)}$,
\begin{equation}\label{eq:erreur}
\begin{aligned}
\|{\mathcal L}(\tau_S-t,X_{\tau_S} |  \tau_S
  > t) - {\mathcal L}(\tau_S,X_{\tau_S} | X_0 \sim \nu_S)\|_{TV} \le C {\rm e}^{-(\lambda_2-\lambda_1) t }
\end{aligned}
\end{equation}
where
$$
\begin{aligned}
&\|{\mathcal L}(\tau_S-t,X_{\tau_S} |  \tau_S
  > t) - {\mathcal L}(\tau_S,X_{\tau_S} | X_0 \sim \nu_S)\|_{TV}\\
&=
\sup_{f,\, \|f\|_{L^\infty} \le 1}  \left| \E (f(\tau_S-t,X_{\tau_S}) | \tau_S
  > t) -
  \E^{\nu_S} (f(\tau_S,X_{\tau_S}))  \right|
\end{aligned}
$$
denotes the total variation norm of the difference between the law of
$(\tau_S-t,X_{\tau_S})$ conditioned to $\tau_S > t$ (for any initial
condition $X_0 \in S$), and the law of
$(\tau_S,X_{\tau_S})$ when $X_0$ is distributed according to $\nu_S$.
The supremum is taken over all bounded functions $f:\R_+
\times \partial S \to \R$, with $L^\infty$-norm smaller than one.

This gives a way to quantify the local equilibration time mentioned in the
introduction of Section~\ref{sec:QSD}, which is the typical time to
get the convergence in~\eqref{eq:QSD_1}: it is of order
$1/(\lambda_2-\lambda_1)$. Of course, this is not a very practical
result since computing the eigenvalues $\lambda_1$ and $\lambda_2$ is in general
impossible. We will discuss in Section~\ref{sec:ParRep_gen} a practical way
to estimate this time.

As a consequence, this result also gives us a way to define a
metastable state: the local equilibration time is of order
$1/(\lambda_2-\lambda_1)$, the exit time is of order $1/\lambda_1$ and
thus, the state $S$ is metastable if
\begin{equation}\label{eq:S_metastable}
\frac{1}{\lambda_1} \gg \frac{1}{\lambda_2-\lambda_1}.
\end{equation}

\subsection{A first discussion on QSD and jump Markov model}\label{sec:QSD_Markov}

Let us now go back to our discussion on the link between the
overdamped Langevin dynamics~\eqref{eq:OL} and the jump Markov
dynamics~\eqref{eq:Z}. Using the first property~\ref{sec:QSD1}, if the
process remains in $S$ for a long time, then it is approximately
distributed according to the QSD, and the error can be quantified
thanks to~\eqref{eq:erreur}. Therefore, to study the exit from
$S$, it is relevant to consider a process starting from the QSD
$\nu_S$ in $S$. Then, the third property~\ref{sec:QSD3} shows that the
exit event can indeed be identified with one step of a Markov jump
process since $\tau_S$ is exponentially distributed and independent of
$X_{\tau_S}$, which are the basic requirements of a move of a Markov jump process
(see Section~\ref{sec:Z}). 

In other words, the QSD $\nu_S$ is the natural initial distribution to choose
in a metastable state $S$ in order to parametrize an underlying jump
Markov model.

In order to be more precise, let us assume that the state $S$ is
surrounded by $I$ neighboring states. The boundary $\partial S$ is
then divided into $I$ disjoint subsets $(\partial S_i)_{i=1, \ldots, I}$, each of 
them associated with
 an exit towards one of the neighboring
states, which we assume to be numbered by $1, \ldots , I$ without loss
of generality: ${\mathcal N}_0=\{1, \ldots, I\}$ (see Figure~\ref{fig:S} for a situation where $I=4$). The
exit event from $S$ is characterized by the pair
 $(\tau_S,\mathcal I)$,
 where $\mathcal I$ is a random variable which gives the next
visited state:
\[\text{for }i=1,\ldots,I,\quad \{\mathcal I=i\}=\{X_{\tau_S}
\in \partial S_i\}.\]
Notice that $\tau_S$ and $\mathcal I$ are by construction independent random variables.
The jump Markov model is
then parametrized as follows. Introduce (see Equation~\eqref{eq:QSD_3}
for the exit point distribution) 
\begin{equation}\label{eq:p}
p(i)= \P(X_{\tau_S} \in \partial S_i) =
-\frac{\displaystyle \int_{\partial S_i} \partial_n u_1 \, d \sigma}{\displaystyle \beta \lambda_1
  \int_S u_1(x)  \, d x},\quad
\text{for } i=1,\ldots,I.
\end{equation}
For each exit region $\partial S_i$, let us define the corresponding rate
\begin{equation}\label{eq:k}
\text{for } i=1,\ldots,I, \, k_{0,i}=\lambda_1 p(i).
\end{equation}
Now, one can check that 
\begin{itemize}
\item The exit time $\tau_S$ is exponentially distributed with
parameter $\sum_{i \in {\mathcal N}_0} k_{0,i}$, in accordance
with~\eqref{eq:T_n}.
\item The next visited state is ${\mathcal I}$,
independent of $\tau_S$ and with law: for $j \in {\mathcal N}_0,  \, \P({\mathcal I}=j)=\frac{k_{0,j}}{\sum_{i \in {\mathcal N}_0} k_{0,i}}$, in
accordance with~\eqref{eq:Y_n}.
\end{itemize}
Let us emphasize again that $\tau_S$ and
$X_{\tau_S}$ are independent random variables, which is a crucial
property to recover the Markov jump model (in~\eqref{eq:T_n}--\eqref{eq:Y_n}, conditionally
on $Y_n$, $T_n$ and $Y_{n+1}$ are indeed independent).

The rates given by~\eqref{eq:k} are exact, in the sense that starting
from the QSD, the law of the exit event from $S$ is exact using this definition
for the transitions to neighboring states. In Section~\ref{sec:EK}, we
will discuss the error introduced when approximating these rates by
the Eyring-Kramers formula~\eqref{eq:EK}.

As a comment on the way we define the rates, let us mention that in
the original works by Kramers~\cite{kramers-40} (see also~\cite{naeh-klosek-matkowsky-schuss-90}), the idea is to
introduce the stationary Fokker-Planck equation with
zero boundary condition (sinks on the boundary of $S$) and with a
source term within $S$ (source in $S$), and to look at the steady
state outgoing current on the boundary $\partial S$. When the process leaves $S$, it is reintroduced in $S$
according to the source term. In general, the stationary state depends
on the source term of course. The difference with the QSD
approach (see~\eqref{eq:QSD_2}) is that we consider the first
eigenvalue of the Fokker-Planck operator. This corresponds to the
following: when the process leaves
$S$, it is reintroduced in $S$ according to the empirical law along
the path of the process in $S$. The interest of this point of view is
that the exit time distribution is exactly exponential (and not
approximately exponential in some small temperature or high barrier regime).

\subsection{Concluding remarks}

The interest of the QSD approach is that it is very general and versatile. The QSD
can be defined for any stochastic process: reversible or
non-reversible, with values in a discrete or a continous state space,
etc, see~\cite{collet-martinez-san-martin-13}. Then, the properties that the exit time is exponentially
distributed and independent of the exit point are satisfied in these very
general situations. 

Let us emphasize in particular that in the
framework of the two dynamics~\eqref{eq:L} and~\eqref{eq:OL} we
consider here, the QSD gives a natural way to define rates to leave a
metastable state, without any small temperature assumption. Moreover,
the metastability may be related to either energetic
barriers or entropic barriers (see in particular~\cite{binder-simpson-lelievre-15} for
numerical experiments in purely entropic cases). Roughly speaking, energetic
barriers correspond to a situation where it is difficult to leave $S$
because it corresponds to the basin of attraction of a local minimum
of $V$ for the gradient dynamics $\dot{x}=-\nabla V(x)$: the process
has to go over an energetic hurdle (namely a saddle point of $V$) to
leave $S$. Entropic barriers are different. They appear when it takes
time for the process to leave $S$ because the exit doors from $S$ are very
narrow. The potential within $S$ may be constant in this case. In practice, entropic barriers are related to steric
constraints in the atomic system. The extreme case for an entropic
barrier is a Brownian motion ($V=0$)
reflected on $\partial S \setminus \Gamma$, $\Gamma \subset \partial
S$ being the small subset of $\partial S$ through which the process
can escape from $S$. For applications in biology for
example, being able to handle both energetic and entropic barriers is
important. 

Let us note that the QSD in $S$ is in general different from the
Boltzmann distribution restricted to $S$: the QSD is zero on the
boundary of $\partial S$ while this is not the case for the Boltzmann distribution.

The remaining of the article is organized as follows. In
Section~\ref{sec:Voter}, we review recent results which show how the QSD can be used to justify
and analyze accelerated dynamics algorithms, and in particular the
parallel replica algorithm. These techniques aim at efficiently
sample the state-to-state dynamics associated with the microscopic
models~\eqref{eq:L} and~\eqref{eq:OL}, using the underlying jump
Markov model to accelerate the sampling of the exit event from
metastable states. In Section~\ref{sec:EK}, we
present new results concerning the justification of the Eyring-Kramers
formula~\eqref{eq:EK} for parametrizing a jump Markov model. The two following
sections are essentially independent of each other and can be read separately.

\section{Numerical aspects: accelerated dynamics}\label{sec:Voter}

As explained in the introduction, it is possible to use the underlying
Markov jump process as a support to accelerate molecular
dynamics. This is the principle of the accelerated dynamics methods
introduced by A.F. Voter in the late
nineties~\cite{voter-97,voter-98,sorensen-voter-00}. These techniques
aim at efficiently simulate the exit event from a metastable state.

Three ideas have
been explored. In the parallel replica
algorithm~\cite{voter-98,perez-uberuaga-voter-15}, the idea is to use
the jump Markov model in order to parallelize the sampling of the exit
event. The principle of the hyperdynamics algorithm~\cite{voter-97} is
to raise the potential within the metastable states in order to
accelerate the exit event, while being able to recover the correct
exit time and exit point distributions. Finally, the temperature
accelerated dynamics~\cite{sorensen-voter-00} consists in simulating
exit events at high temperature, and to extrapolate them at low
temperature using the Eyring-Kramers law~\eqref{eq:EK}. In this paper,
for the sake of conciseness, we concentrate on the analysis of the
parallel replica method, and we refer to the
papers~\cite{lelievre-nier-15,aristoff-lelievre-14} for an
analysis of hyperdynamics and temperature accelerated dynamics. See
also the recent review~\cite{lelievre-15} for a detailed presentation.

\subsection{The parallel replica method}\label{sec:ParRep}

In order to present the parallel replica method, we need to introduce
a partition of the configuration space $\R^d$ to describe the
states. Let us denote by
\begin{equation}\label{eq:S}
{\mathcal S}: \R^d \to \N
\end{equation}
a function which associates to a configuration $x \in \R^d$ a state
number ${\mathcal S}(x)$. We will discuss below how to choose in
practice this function ${\mathcal S}$.
The aim of the parallel replica method (and
actually also of hyperdynamics and temperature accelerated dynamics)
is to generate very efficiently a trajectory $(S_t)_{t \ge 0}$ with
values in $\N$ which has approximately the same law as the
state-to-state dynamics $({\mathcal
  S}(X_t))_{t \ge 0}$ where $(X_t)_{t \ge 0}$
follows~\eqref{eq:OL}. The states are the level sets of ${\mathcal
  S}$. Of course, in general, $({\mathcal
  S}(X_t))_{t \ge 0}$ is not a Markov process, but it is close to
Markovian if the level sets of ${\mathcal S}$ are metastable regions, see Sections 2.2
and 2.5. The idea is to check and then use metastability of the states
in order to efficiently generate the exit events.

As explained above, we present for the sake of simplicity the algorithm in
the setting of the overdamped Langevin dynamics~\eqref{eq:OL}, but the
algorithm and the discussion below can be generalized to the Langevin
dynamics~\eqref{eq:L}, and actually to any Markov dynamics, as soon as
a QSD can be defined in each state.

The parallel replica algorithm consists in iterating three steps:
\begin{itemize}
\item {\em The decorrelation step}: In this step, a reference replica
  evolves according to the original dynamics~\eqref{eq:OL}, until it
  remains trapped for a time $t_{corr}$ in one of the states
  ${\mathcal S}^{-1}(\{n\})=\{x \in
  \R^d, \, {\mathcal S}(x)=n\}$, for $n\in \N$. The parameter
  $t_{corr}$ should be chosen by the user, and may depend on the
  state. During this step, no error is made, since the reference
  replica evolves following the original
  dynamics (and there is of course no computational gain compared to a
  naive direct numerical simulation).  Once the reference replica has been
  trapped in one of the states (that we denote generically by $S$ in
  the following two steps) for a time $t_{corr}$, the aim is to generate very
  efficiently the exit event. This is done in two steps.
\item {\em The dephasing step}: In this preparation step, $(N-1)$
  configurations are generated within $S$ (in addition to the one
  obtained form the reference replica) as follows. Starting from the
  position of the reference replica at the end of the decorrelation
  step, some trajectories are simulated in parallel for a time
  $t_{corr}$. For each trajectory, if it remains within $S$ over the
  time interval of length $t_{corr}$, then its end point is
  stored. Otherwise, the trajectory is discarded, and a new attempt to get
  a trajectory remaining in $S$ for a time $t_{corr}$ is made. This
  step is pure overhead. The objective is only to get $N$
  configurations in $S$ which will be used as initial conditions in
  the parallel step.
\item {\em The parallel step}: In the parallel step, $N$ replicas are
  evolved independently and in parallel, starting from the initial conditions generated
  in the dephasing step, following the original dynamics~\eqref{eq:OL}
  (with independent driving Brownian motions). This step ends as soon
  as one of the replica leaves $S$. Then, the simulation clock is
  updated by setting the residence time in the state $S$ to $N$ (the
  number of replicas) times the exit time of the first replica which
  left $S$. This replica now becomes the reference replica, and one
  goes back to the decorrelation step above.
\end{itemize}
The computational gain of this algorithm is in the parallel step,
which (as explained below) simulates the exit event in a wall clock
time $N$ times smaller  in average than what would have been necessary to see the
reference walker leaving $S$. This of course requires a parallel
architecture able to handle $N$ jobs in parallel\footnote{For a
  discussion on the parallel efficiency, communication and
  synchronization, we refer to the papers~\cite{voter-98,perez-uberuaga-voter-15,le-bris-lelievre-luskin-perez-12,binder-simpson-lelievre-15}.}. This algorithm can be seen as a way to
parallelize in time the simulation of the exit event, which is not trivial because of the sequential nature of time.

Before we present the mathematical analysis of this method, let us
make a few comments on the choice of the function ${\mathcal
  S}$. In the original papers~\cite{voter-98,perez-uberuaga-voter-15},
the idea is to define states as the basins of attraction of the local
minima of $V$ for the gradient dynamics $\dot{x}=-\nabla V(x)$. In this context, it is
important to notice that the states do not need to be defined {\em a
  priori}: they are numbered as the process evolves and discovers new
regions (namely new local minima of $V$ reached by the gradient
descent). This way to define ${\mathcal S}$ is well suited for
applications in material sciences, where barriers are essentially
energetic barriers, and the local minima of $V$ indeed correspond to
different macroscopic states. In other applications, for example in
biology, there may be too many local minima, not all of them being
significant in terms of macroscopic states. In that case, one could
think of using a few degrees of freedom (reaction coordinates) to
define the states, see for
example~\cite{kum-dickson-stuart-uberuaga-voter-04}. Actually, in the
original work by Kramers~\cite{kramers-40}, the states are also
defined using reaction coordinates, see the discussion in~\cite{hanggi-talkner-barkovec-90}. The important
outcome of the mathematical analysis below is that, whatever the
choice of the states, if one is able to define a correct correlation
time $t_{corr}$ attached to the states, then the algorithm is consistent. We will discuss in Section~\ref{sec:ParRep_math}
how large $t_{corr}$ should be theoretically, and in
Section~\ref{sec:ParRep_gen} how to estimate it in practice.

Another important remark is that one actually does not need a
partition of the configuration space to apply this algorithm. Indeed,
the algorithm can be seen as an efficient way to simulate the exit event
from a metastable state $S$. Therefore, the algorithm could be applied even if no partition of the state space is available, but only an ensemble of disjoint subsets of the configuration space. The algorithms could then be used to simulate efficiently exit events from these states, if the system happens to be trapped in one of them.

\subsection{Mathematical analysis}\label{sec:ParRep_math}

Let us now analyze the parallel replica algorithm described above,
using the notion of quasi-stationary distribution. In view of the
first property~\ref{sec:QSD1} of the QSD, the decorrelation step is
simply a way to decide wether or not the reference replica remains
sufficiently long in one of the states so that it can be considered as
being distributed according to the QSD. In view of~\eqref{eq:erreur},
the error is of the order of $\exp(-(\lambda_2-\lambda_1) \, t_{corr})$
so that
$t_{corr}$ should be chosen of the order of $1/(\lambda_2 -
\lambda_1)$ in order for the exit event of the reference walker which
remains in $S$ for a time $t_{corr}$ to be statistically close to the
exit event generated starting from the QSD. 

Using the same arguments, the dephasing step is nothing but a rejection algorithm to generate
many configurations in $S$ independently and identically distributed
with law the QSD $\nu_S$ in $S$. Again, the distance to the QSD of the
generated samples can be
quantified using~\eqref{eq:erreur}.

Finally, the parallel step generates an exit event which is exactly
the one that would have been obtained considering only one
replica. Indeed, up to the error quantified in~\eqref{eq:erreur}, all
the replica are i.i.d. with initial condition the QSD
$\nu_S$. Therefore, according to the third property~\ref{sec:QSD3} of
the QSD, their exit times $(\tau^n_S)_{n \in \{1, \ldots N\}}$ are i.i.d. with law an
exponential distribution with parameter $\lambda_1$ ($\tau^n_S$ being the exit time
of the $n$-th replica) so that
\begin{equation}\label{eq:Nexp}
N \min_{n \in \{1, \ldots, N\}} (\tau^n_S) \stackrel{\mathcal L}{=}
\tau^1_S.
\end{equation}
This explains why the exit time of the first replica which leaves $S$
needs to be multiplied by the number of replicas $N$. This also shows
why the parallel step gives a computational gain in terms of wall
clock: the time required to simulate the exit event is divided by $N$
compared to a direct numerical simulation.
Moreover, since starting from the QSD, the exit time and the exit
point are independent, we also have
$$X^{I_0}_{\tau^{I_0}_S} \stackrel{\mathcal L}{=}
    X^1_{\tau^1_S},$$
where $(X^n_t)_{t \ge 0}$ is the $n$-th replica and $I_0=\arg\min_{n \in \{1, \ldots, N\}}(\tau^n_S)$ is the index
of the first replica which exits $S$. The exit
point of the first replica which exits $S$ is statistically the same as
the exit point of the reference walker. Finally, by the independence property of
exit time and exit point, one can actually combine the two former
results in a single equality in law on couples of random variables, which shows that the parallel step is
statistically exact:
$$\left(N \min_{n \in \{1, \ldots, N\}} (\tau^n_S) ,
X^{I_0}_{\tau^{I_0}_S}\right)
\stackrel{\mathcal L}{=} (\tau^1_S, X^1_{\tau^1_S}).$$
As a remark, let us notice that in practice, discrete-time processes
are used (since the Langevin or overdamped Langevin dynamics are discretized in time). Then, the exit
times are not exponentially but geometrically distributed. It is
however possible to generalize the formula~\eqref{eq:Nexp} to this
setting by using the following fact:  if $(\sigma_n)_{n \in \{1, \ldots
  N\}}$ are i.i.d. with geometric law, then
$N \left( \min(\sigma_1,
\ldots, \sigma_N)-1 \right) + \min \left( n \in \{1, \ldots ,N \} , \, \sigma_n =  \min(\sigma_1,
\ldots, \sigma_N) \right) \stackrel{{\mathcal L}}{=}\sigma_1$. We refer
to~\cite{aristoff-lelievre-simpson-14} for more details. 

This analysis shows that the parallel replica is a very versatile algorithm. In
particular it applies to both energetic and entropic barriers, and
does not assume a small temperature regime (in contrast with the
analysis we will perform in Section~\ref{sec:EK}). The
only errors introduced in the algorithm are related to the rate of
convergence to the QSD of the process conditioned to stay in the
state. The algorithm will be efficient if the convergence time to the
QSD is small compared to the exit time (in other words, if the states
are metastable). Formula~\eqref{eq:erreur}
gives a way to quantify the error introduced by the whole
algorithm. In the limit $t_{corr} \to \infty$, the algorithm generates
exactly the correct exit event. However,~\eqref{eq:erreur} is not very
useful to choose $t_{corr}$ in practice since it is not
possible to get accurate estimates of $\lambda_1$ and $\lambda_2$ in
general. We will present in the next section a practical way to
estimate $t_{corr}$. 

Let us emphasize that this analysis gives some error bound on the
accuracy of the {\em state-to-state dynamics} generated by the
parallel replica algorithm, and not only on the invariant measure, or
the evolution of the time marginals.






\subsection{Recent developments on the parallel replica algorithm}\label{sec:ParRep_gen}

In view of the previous mathematical analysis, an important practical
question is how to choose the correlation time $t_{corr}$. In the original
papers~\cite{voter-98,perez-uberuaga-voter-15}, the correlation time
is estimated assuming that an harmonic approximation is accurate. In~\cite{binder-simpson-lelievre-15}, we propose
another approach which could be applied in more general settings. The
idea is to use two ingredients:
\begin{itemize}
\item The Fleming-Viot particle process~\cite{ferrari-maric-07}, which
  consists in
  $N$ replicas $(X^1_t, \ldots,X^N_t)_{t \ge 0}$ which are evolving and
  interacting in such a way that
  the empirical distribution $\frac{1}{N} \sum_{n=1}^N \delta_{X^n_t}$
  is close (in the large $N$ limit) to the law of the process $X_t$
  conditioned on $t<\tau_S$.
\item The Gelman-Rubin convergence diagnostic~\cite{gelman-rubin-92} to estimate the
  correlation time as the convergence time to a stationary state for
  the Fleming-Viot particle process.
\end{itemize}
Roughly speaking, the Fleming-Viot process consists in following the
original dynamics~\eqref{eq:OL} independently for each replica, and, each time one of the replicas leaves the domain $S$, another one
taken at random is duplicated. The Gelman-Rubin convergence diagnostic
consists in comparing the average of a given observable over replicas
at a given time, with the average of this observable over time and
replicas: when the two averages are close (up to a tolerance, and for
well chosen observables), the process is considered at stationarity.

Then, the generalized parallel replica algorithm introduced
in~\cite{binder-simpson-lelievre-15} is a modification of the original
algorithm where, each time the reference replica enters a new state, a
Fleming-Viot particle process is launched using $(N-1)$ replicas simulated in
parallel. Then the decorrelation step consists in the following: if
the reference replica leaves $S$ before the Fleming-Viot particle
process reaches stationarity, then a new decorrelation step starts
(and the replicas generated by the Fleming-Viot particle are discarded);
if otherwise the Fleming-Viot particle
process reaches stationarity before the reference replica leaves $S$,
then one proceeds to the parallel step. Notice indeed that the final positions
of the replicas simulated by the Fleming-Viot particle process can be
used as initial conditions for the processes in the parallel
step. This procedure thus avoids the choice of a $t_{corr}$ {\em a
  priori}: it is in some sense estimated on the fly. For
more details, discussions on the correlations included by the
Fleming-Viot process between the replicas, and numerical experiments
(in particular in cases with purely entropic barriers), we refer to~\cite{binder-simpson-lelievre-15}.

\subsection{Concluding remarks}\label{sec:ParRep_conc}

We presented the results in the context of the overdamped Langevin
dynamics~\eqref{eq:OL}, but the algorithms straightforwardly apply to
any stochastic Markov dynamics as soon as a QSD exists (for example
Langevin dynamics for a bounded domain, see~\cite{nier-14}).

The QSD approach is also useful to analyze the two other accelerated
dynamics: hyperdynamics~\cite{lelievre-nier-15} and temperature accelerated
dynamics~\cite{aristoff-lelievre-14}. Typically, one expects better
speed up with these algorithms than with parallel replica, but at the
expense of larger errors and more stringent assumptions (typically
energetic barriers, and small temperature regime),
see~\cite{lelievre-15} for a review paper.
Let us mention in particular that the mathematical analysis of the
temperature accelerated dynamics algorithms requires to prove that the
distribution for the next visited state predicted using the
Eyring-Kramers formula~\eqref{eq:EK} is correct, as explained
in~\cite{aristoff-lelievre-14}. The next section is thus also
motivated by the development of an error analysis for temperature accelerated dynamics.

Let us finally mention that in these algorithms, the way to relate the
original dynamics to a jump Markov process is by looking at
$({\mathcal S}(X_t))_{t \ge 0}$ (or $({\mathcal S}(q_t))_{t \ge 0}$
for~\eqref{eq:L}). As already mentioned, this is not a Markov
process, but it is close to Markovian if the level sets of ${\mathcal
  S}$ are metastable sets, see Sections~\ref{sec:QSD1} and~\ref{sec:erreur}. In particular, in the
parallel replica algorithm above, the non-Markovian effects (and in particular
the recrossing at
the boundary between two states) are taken into account using the
decorrelation step, where the exact process is used in these
intermediate regimes between long sojourns in metastable states. As
already mentioned above (see the discussion on the map ${\mathcal S}$
at the end of Section~\ref{sec:ParRep}), another idea is to  introduce an
ensemble of disjoint subsets $(M_i)_{i \ge 0}$ and to project the
dynamics $(X_t)_{t \ge 0}$ (or $(q_t)_{t \ge 0}$) onto a discrete state-space dynamics by considering {\em the last
visited
milestone} \cite{buchete-hummer-08,schuette-noe-lu-sarich-vanden-einjden-11}. Notice
that these subsets do not create a partition of the state space. They are sometimes called
milestones \cite{faradjian-elber-04}, target sets or core sets \cite{schuette-noe-lu-sarich-vanden-einjden-11} in the literature. The natural parametrization of the
underlying jump process is then to consider, starting from a
milestone (say $M_0$), the time to reach any of the other milestones
($(M_j)_{j \neq 0}$) and the index of the next visited
milestone. This requires us 
to study the reactive paths among the
milestones, for which many techniques have been
developed, as already presented in the introduction. Let us now discuss the Markovianity of the projected
dynamics. On the one hand, in the limit of very small
milestones\footnote{One could think of 
one-dimensional overdamped Langevin
  dynamics, with
  milestones defined as points: in this case the sequence of visited
  points is Markovian.}, the sequence of visited states
({\em i.e.} the skeleton of the projected process) is naturally Markovian (even
though the transition time is not necessarily exponential), but the
description of the underlying continuous state space dynamics is very
poor (since the information of the last visited milestone is not very
informative about
 the actual state of the system). On the other hand, taking
larger milestones, the projected process is close to a Markov process
under some metastability assumptions with respect to these
milestones. We refer to \cite{sarich-noe-schuette-10,schuette-sarich-13,bovier-den-hollander-15} for a mathematical analysis.

\section{Theoretical aspects: transition state theory and
  Eyring-Kramers formula}\label{sec:EK}

In this section, we explore some theoretical
counterparts of the QSD approach to study metastable stochastic
processes. We concentrate on the overdamped Langevin
dynamics~\eqref{eq:L}. The generalization of the mathematical approach
presented below to the Langevin dynamics would require some extra
work. 

We would like to justify the procedure
described in the introduction to build jump Markov models, and which
consists in (see for example~\cite{cameron-2014b,wales-03,voter-05}): (i) looking for all local minima and saddle points
separating the local minima of the function $V$;
(ii) connecting two minima which can be linked by a path going through a
single saddle point, and parametrizing a
jump between these two minima using the rate given by the Eyring-Kramers
formula~\eqref{eq:EK}. 
More precisely, we concentrate on the accuracy of the sampling of the
exit event from a metastable state using the jump Markov model. The
questions we ask are the following: if a set $S$
containing a single local minimum of $V$ is metastable for the dynamics~\eqref{eq:OL} (see the discussion in
Section~\ref{sec:QSD1} and formula~\eqref{eq:S_metastable}), is the
exit event predicted by the jump Markov model built using the
Eyring-Kramers formula correct? What is the error induced by this
approximation?

As already explained in Section~\ref{sec:QSD_Markov}, if $S$ is
metastable, one can assume that the stochastic process $(X_t)_{t \ge
  0}$ satisfying~\eqref{eq:OL} starts under the QSD $\nu_S$ (the error
being quantified by~\eqref{eq:erreur}) and then,
the exit time is exponentially distributed and independent of the exit point.
Thus, two fundamental properties of the jump Markov model are
satisfied. It only remains to prove that the rates associated
with the exit event for $(X_t)_{t \ge 0}$ (see formula~\eqref{eq:k})
can be accurately approximated by the Eyring-Kramers
formulas~\eqref{eq:EK}. As will become clear below, the analysis holds for energetic barriers in the small temperature regime
$\beta \to \infty$.

In this section, we only sketch the proofs of our results, which are
quite technical. For a more
detailed presentation, we refer to~\cite{di-gesu-le-peutrec-lelievre-nectoux-16}.

\subsection{A review of the literature}\label{sec:EK_litt}

Before presenting our approach, let us discuss the mathematical
results in the literature aiming at justifying the
Eyring-Kramers rates. See also the review article~\cite{berglund-13}.

Some authors adopt a global approach: they look at the spectrum
associated with the infinitesimal generator of the dynamics on the whole
configuration space, and
they compute the small eigenvalues in the small temperature regime
$\beta \to \infty$. It can be shown that there are exactly $m$ small
eigenvalues, $m$ being the number of local minima of $V$, and that these
eigenvalues satisfy the Eyring-Kramers law~\eqref{eq:EK}, with an
energy barrier $V(z_k)-V(x_k)$. Here, the
saddle point $z_k$ attached to the local minimum $x_k$ is defined by\footnote{It is here
  implicitly assumed that the inf sup value is attained at a single
  saddle point~$z_k$.}
$$V(z_k)=\inf_{\gamma \in {\mathcal P}(x_i,B_i)} \sup_{t \in [0,1]}
V(\gamma(t))$$
where ${\mathcal P}(x_i,B_i)$ denotes the set of continuous paths from
$[0,1]$ to $\R^d$ such that $\gamma(0)=x_i$ and $\gamma(1) \in B_i$ with
$B_i$ the union of small balls around local minima lower in energy
than $x_i$. For the dynamics~\eqref{eq:OL}, we refer for example to the
work~\cite{helffer-klein-nier-04} based on semi-classical analysis
results for Witten Laplacian and the
articles~\cite{bovier-eckhoff-gayrard-klein-04,bovier-gayrard-klein-05,eckhoff-05}
where a potential theoretic approach is adopted. In the latter
results, a connexion is made between the small eigenvalues and mean
transition times between metastable states. Let us also mention the
earlier results~\cite{miclo-95,holley-kusuoka-stroock-89}.
For the dynamics~\eqref{eq:L}, similar results are obtained
in~\cite{herau-hitrick-sjostrand-11}. These spectral approaches give the
cascade of relevant time scales to reach from a local minimum any
other local minimum which is lower in energy. They do not give any information about the typical time scale to go from one local minimum to any other local minimum (say from the global minimum to the second lower minimum). 
These global approaches can be used to build jump Markov models using a
Galerkin projection of the infinitesimal generator onto the first $m$
eigenmodes, which gives an excellent approximation of the
infinitesimal generator. This has been extensively investigated by
Sch\"utte\footnote{In fact, Sch\"utte et al. look at the eigenvalues
  close to 1 for the so-called transfer operator $P_t = {\rm e}^{tL}$
  (for a well  chosen lag time $t > 0$), which is equivalent to looking at the small positive eigenvalues of $-L$} and his collaborators~\cite{schuette-sarich-13}, starting with the seminal work~\cite{schuette-98}.

In this work, we are interested in a local approach, namely in the
study of the exit event from a given metastable state $S$. In this
framework, the most famous approach to analyze the exit event is the
large deviation theory~\cite{freidlin-wentzell-84}. In the small
temperature regime, large deviation results
provide the exponential rates~\eqref{eq:EK}, but without the
prefactors and without error bounds. It can also be proven that the exit time is exponentially
distributed in this regime, see~\cite{day-83}. For the
dynamics~\eqref{eq:OL}, a typical result on the exit point distribution
is the following (see~\cite[Theorem 5.1]{freidlin-wentzell-84}): for all $S'
\subset\subset S$, for any $\gamma >0$, for any $\delta > 0$, there
exists $\delta_0 \in (0,\delta]$ and $\beta_0 > 0$ such that for
all $\beta \ge \beta_0$, for all $x \in S'$ and for all $y
\in \partial S$,
\begin{equation}\label{eq:LD}
{\rm e}^{-\beta ( V(y)-V(z_1) + \gamma)}\le \P^{x} (X_{\tau_S} \in
{\mathcal V}_{\delta_0}(y))  \le
{\rm e}^{-\beta ( V(y)-V(z_1)- \gamma)}
\end{equation}
where ${\mathcal V}_{\delta_0}(y)$ is a $\delta_0$-neighborhood of $y$ in
$\partial S$. Besides, let us also mention formal approaches to study the exit time and the exit point
distribution that have been proposed by Matkowsky, Schuss and collaborators
in~\cite{matkowsky-schuss-77,naeh-klosek-matkowsky-schuss-90,schuss-09}
and by Maier and Stein in~\cite{maier-stein-93}, using formal
expansions for singularly perturbed elliptic equations. Some of the
results cited above actually consider more general dynamics
than~\eqref{eq:OL} (including~\eqref{eq:L}), see
also~\cite{,bouchet-reygner-15} for a recent contribution in that
direction. One of the interests of the large deviation
approach is actually to be sufficiently robust to apply to rather general dynamics.

Finally, some authors prove the convergence to a jump Markov process using a
rescaling in time. See for example~\cite{kipnis-newman-85} for
a one-dimensional diffusion in a double well, and~\cite{galves-olivieri-vares-87,mathieu-95} for a
similar problem in larger dimension. In \cite{sugiura-95}, a rescaled
in time diffusion process converges to a jump Markov process living on
the global minima of the potential $V$, assuming they are separated
by saddle points having the same heights.

There are thus many mathematical approaches to derive the Eyring-Kramers
formula. In particular, a lot of works are devoted to the computation
of the rate between
two metastable states, but very few discuss the use
of the combination of these rates to build a jump Markov model between
metastable states.
To the best of our knowledge, none of
these works quantify rigorously the error introduced by the use of the Eyring-Kramers formulas and a jump Markov
process to model the transition from one state to all the neighboring
states. Our aim in this section is to present such a mathematical
analysis, using local versions of the spectral approaches mentioned
above.
Our approach is local, justifies the Eyring-Kramers formula with the
prefactors and provides error estimates. It uses techniques developed in particular in the previous works~\cite{helffer-nier-06,le-peutrec-10}.
These results generalize the
results in dimension~1 in Section 4 of~\cite{aristoff-lelievre-14}.

\subsection{Mathematical result}

Let us consider the dynamics~\eqref{eq:OL} with an initial condition
distributed according to the QSD $\nu_S$ in a domain $S$. We assume
the following:
\begin{itemize}
\item The domain $S$ is an open smooth bounded domain in $\R^d$.
\item The function $V:\overline{S} \to \R$ is a Morse function with a single critical
point $x_1$. Moreover, $x_1 \in S$ and $V(x_1)=\min_{\overline{S}} V$. 
\item The normal derivative $\partial_n V$ is strictly positive on $\partial
  S$, and $V|_{\partial S}$ is a Morse function with local minima
  reached at $z_1, \ldots , z_I$ with $V(z_1) < V(z_2) < \ldots <
  V(z_I)$.
\item The height of the barrier is large compared to the saddle points
  heights discrepancies: $V(z_1)-V(x_1) > V(z_I)-V(z_1)$.
\item For all $i \in \{1, \ldots I\}$, consider $B_{z_i}
  \subset \partial S$ the basin of
  attraction for the dynamics in the boundary $\partial S$: $\dot{x} = -\nabla_T V(x)$ (where
  $\nabla_T V$ denotes the tangential gradient of $V$ along the
  boundary $\partial S$). Assume that
\begin{equation}\label{eq:agmon}
\inf_{z \in B_{z_i}^c} d_a (z,z_i) > V(z_I)-V(z_1)
\end{equation}
where $B_{z_i}^c=\partial S \setminus B_{z_i}$.
\end{itemize}
Here, $d_a$ is the Agmon distance:
$$d_a(x,y) = \inf_{\gamma \in \Gamma_{x,y}} \int_0^1 g(\gamma(t)) |\gamma'(t)| \, dt$$
where $g=\left\{ \begin{aligned} &|\nabla V| \text { in $S$} \\ &|\nabla_T
    V|  \text{ in $\partial S$} \end{aligned} \right.$, and the
infimum is over the set $\Gamma_{x,y}$ of all piecewise $C^1$ paths
$\gamma:[0,1] \to \overline{S}$ such that $\gamma(0)=x$ and
$\gamma(1)=y$.
The Agmon distance is useful in order to measure the decay of
eigenfunctions away from critical points. These are the so-called
semi-classical Agmon estimates,
see~\cite{simon-84,helffer-sjostrand-84}.

Then,  in the limit $\beta \to \infty$,
the exit rate is (see also~\cite{helffer-nier-06}) $$\lambda_1= \sqrt{\frac{\beta}{2\pi}}
\partial_nV(z_1) \frac{  \displaystyle\sqrt{ \det  (\nabla^2V)   (x_1) }
}{    \displaystyle\sqrt{ \det (\nabla^2 V_{|\partial S})   (z_1) }
}{\rm e}^{-\beta(V(z_1)-V(x_1))}( 1+ O(\beta^{-1}) ).$$
Moreover, for any open set $\Sigma_i $ containing $z_i$ such that $\overline\Sigma_i \subset B_{z_i}$, 
\begin{equation}\label{eq:dnu1}
\frac{\displaystyle \int_{\Sigma_i}   \partial_{n} u_1 \,
  d\sigma}{\displaystyle \int_S u_1(x) \, dx}=- A_i (\beta)
{\rm e}^{-\beta (V(z_i)-V(x_1))}( 1+ O(\beta^{-1}) ),
\end{equation}
 where  $$A_i(\beta)=\frac{\beta^{3/2}}{\sqrt{2\pi}} \partial_n
   V(z_i)  \frac{    \displaystyle \sqrt{\det (\nabla^2 V)   (x_1)} }{  \displaystyle   \sqrt{ {\rm
       det } (\nabla^2 V|_{\partial S})   (z_i) }  }.$$
Therefore,
\begin{align}
 p(i)&=\P^{\nu_S} (X_{\tau_S} \in \Sigma_i) \nonumber \\
&=   
\frac{\partial_n V(z_i) \sqrt{ \det (\nabla^2 V|_{\partial S })
    (z_1) } }{\partial_n V(z_1) \sqrt{ \det (\nabla^2
    V|_{\partial S })   (z_i) }}   {\rm e}^{-\beta(V(z_i)-V(z_1))} (   1+   
 O(\beta^{-1})  ) \label{eq:EK-p}
\end{align}
and (see Equation~\eqref{eq:k} for the definition of the exit rates)
\begin{align}
k_{0,i}&=\lambda_1 p(i) \nonumber \\
&=\widetilde{\nu}^{OL}_{0,i} {\rm e}^{-\beta(V(z_i)-V(x_1))} (   1+   
 O(\beta^{-1})  ) \label{eq:k0i}
\end{align}
where the prefactors $\widetilde{\nu}^{OL}_{0,i}$ are given by
\begin{equation}\label{eq:nu_OL_tilde}
\widetilde{\nu}^{OL}_{0,i}=\sqrt{\frac{\beta}{2\pi}}  \partial_nV(z_i) \frac{ \sqrt{\det (\nabla^2 V)   (x_1) } } {  \sqrt{ \det
     (\nabla^2 V_{|\partial S})   (z_i) }  }.
\end{equation}
We refer to~\cite{di-gesu-le-peutrec-lelievre-nectoux-16} for more
details, and other related results.

As stated in the assumptions, these rates are obtained assuming $\partial_n V >
0$ on $\partial S$: the local minima $z_1, \ldots ,z_I$ of $V$ on
$\partial S$ are therefore
not saddle points of $V$ but so-called {\em generalized saddle points} (see~\cite{helffer-nier-06,le-peutrec-10}). In a future work, we intend to extend these results to the case where
the points $(z_i)_{1 \le i \le I}$ are saddle points of~$V$, in which
case we expect to prove the same result~\eqref{eq:k0i} for the exit
rates, with the prefactor $\widetilde{\nu}^{OL}_{0,i}$ being
$\displaystyle\frac{1}{\pi}|\lambda^-(z_j)| \frac{\displaystyle \sqrt{\det
    (\nabla^2 V)(x_1)}}{\displaystyle\sqrt{|\det (\nabla^2
    V)(z_j)|}}$ (this formula can be obtained using formal expansions on
the exit time and the Laplace's method). Notice that the latter formula differs from~\eqref{eq:nu_OL} by a
multiplicative factor $1/2$ since $\lambda_1$ is the exit rate
from $S$ and not the transition rate to one of the neighboring state (see the
remark on page 408 in~\cite{bovier-eckhoff-gayrard-klein-04} on this multiplicative factor $1/2$ and the results on
asymptotic exit times in~\cite{maier-stein-93} for example). This
factor is due to the fact that once on the saddle point, the process
has a probability one half to go back to $S$, and a probability one half
to effectively leave $S$. This
multiplicative factor does not have any influence on the law of the
next visited state which only involves ratio of the rates $k_{0,i}$,
see Equation~\eqref{eq:Y_n}.

\subsection{Discussion of the result}

As already discussed above, the interest of these results is that they
justify the use of the Eyring-Kramers formula to model the exit
event using a jump Markov model. They give in particular the relative
probability to leave $S$ through each of the local minima $z_i$ of $V$ on the
boundary $\partial S$. Moreover, we obtain an estimate of the
relative error  on the exit probabilities (and
not only on the logarithm of the exit probabilities as
in~\eqref{eq:LD}): it is of order $\beta^{-1}$, see Equation~\eqref{eq:EK-p}.

The importance of obtaining a result including the prefactors in the
rates is illustrated by the following result, which is also proven
in~\cite{di-gesu-le-peutrec-lelievre-nectoux-16}.  Consider a simple
situation with only two local
minima $z_1$ and $z_2$ on the boundary (with as above $V(z_1) < V(z_2)$). Compare the two exit probabilities:
\begin{itemize}
\item The probability to leave
through $\Sigma_2$ such that $\overline{\Sigma_2} \subset B_{z_2}$
and  $z_2 \in \Sigma_2$;
\item The probability to leave through $\Sigma$ such that
  $\overline{\Sigma} \subset B_{z_1}$ and
$\inf_{\Sigma} V=V(z_2)$.
\end{itemize}
By classical results from the large deviation theory (see for
example~\eqref{eq:LD}) the probability to exit through $\Sigma$ and
$\Sigma_2$ both scale like a prefactor times ${\rm
  e}^{-\beta(V(z_2)-V(z_1))}$:  the difference can only be read from the
prefactors. Actually, it can be proven that, in the limit $\beta \to\infty$,
$$\frac{\P^{\nu_S}(X_{\tau_S} \in \Sigma)}{\P^{\nu_S}(X_{\tau_S} \in \Sigma_2)} = O (
\beta^{-1/2} ).
$$
The probability to leave through $\Sigma_2$ (namely through the generalized saddle point $z_2$)
is thus much larger than through $\Sigma$, even though the two regions
are at the same height.
This result explains why the local minima
of $V$ on the boundary (namely the generalized saddle points) play
such an important role when studying the exit event.

\subsection{Sketch of the proof}

In view of the formulas~\eqref{eq:p} and~\eqref{eq:k}, we would like to
identify the asymptotic behavior of the small eigenvalue $\lambda_1$ and of
the normal derivative $\partial_n u_1$ on $\partial S$ in the limit $\beta \to
\infty$. We recall that $(\lambda_1,u_1)$ are defined by the
eigenvalue problem~\eqref{eq:QSD_2}.
In order to work in the classical setting for Witten Laplacians, we
make a unitary transformation of the original eigenvalue problem. Let us consider $v_1=u_1 \exp(\beta V)$, so that 
\begin{equation}\label{eq:L0}
\left\{
\begin{aligned}
L^{(0)} v_1 &= - \lambda_1 v_1 \text{ on $S$,}\\
v_1 & = 0 \text{ on $\partial S$,}
\end{aligned}
\right.
\end{equation}
where $L^{(0)}=\beta^{-1} \Delta - \nabla V \cdot \nabla $
is a self adjoint operator on $L^2(\exp(-\beta V))$. We would like to
study, in the small temperature regime $\partial_n u_1 = \partial_n
v_1 {\rm e}^{-\beta V}$  on $\partial S$ (since $u_1=0$ on $\partial S$). Now, observe that $\nabla v_1$
satisfies
\begin{equation}\label{eq:L1}
\left\{
\begin{aligned}
L^{(1)} \nabla v_1 &= - \lambda_1 \nabla v_1 \text{ on $S$,}\\
\nabla_T v_1 & = 0 \text{ on $\partial S$,}\\
(\beta^{-1} \div - \nabla V \cdot) \nabla v_1 & = 0 \text{ on $\partial S$,}\\
\end{aligned}
\right.
\end{equation}
where
$$L^{(1)}= \beta^{-1} \Delta - \nabla V \cdot \nabla - {\rm Hess}(V)$$
is an operator acting on $1$-forms (namely on vector fields).
Therefore $\nabla v_1$ is an eigenvector (or an eigen-$1$-form)
of  the operator $-L^{(1)}$ with tangential Dirichlet boundary conditions (see~\eqref{eq:L1}), associated with the small eigenvalue $\lambda_1$. It is known (see for
example~\cite{helffer-nier-06}) that in our geometric setting $-L^{(0)}$ admits exactly one
eigenvalue smaller than $\beta^{-1/2}$, namely $\lambda_1$ with
associated eigenfunction $v_1$ (this is because $V$ has only one local
minimum in $S$) and that $-L^{(1)}$ admits exactly $I$
eigenvalues smaller than $\beta^{-1/2}$ (where, we recall, $I$ is the
number of local minima of $V$ on $\partial S$). Actually, all these small
eigenvalues are exponentially small in the regime $\beta \to \infty$,
the larger eigenvalues being bounded from below by a constant in this regime.
The idea is then to construct an appropriate basis (with eigenvectors
localized on the generalized saddle points, see the quasi-modes below) of the eigenspace associated with small eigenvalues for $L^{(1)}$, and then
to decompose $\nabla v_1$ along this basis.

The second step (the most technical one actually) is to build so-called {\em quasi-modes} which approximate the eigenvectors of
$L^{(0)}$ and $L^{(1)}$ associated with small eigenvalues in the
regime $\beta \to \infty$. A good approximation of $v_1$ is actually
simply $\tilde{v}= Z \, \chi_{S'}$ where $S'$ is an open set such
that $\overline{S'} \subset S$,
$\chi_{S'}$ is a smooth function with compact support in $S$ and equal to one on $S'$,  and $Z$ is
a normalization constant such that $\|\tilde v \|_{L^2({\rm e}^{-\beta V})}
= 1$. The difficult part is to
build an approximation of the eigenspace ${\rm Ran}
\left(1_{[0,\beta^{-1/2}]}(-L^{(1)})\right)$, where
$1_{[0,\beta^{-1/2}]}(-L^{(1)})$ denotes the spectral projection of
$(-L^{(1)})$ over eigenvectors associated with eigenvalues in the
interval $[0,\beta^{-1/2}]$. Using auxiliary simpler
eigenvalue problems and WKB expansions around each of the local minima
$(z_i)_{i =1, \ldots ,I}$, we are able to build 
$1$-forms $(\psi_i)_{i=1,\ldots,I}$ such that
${\rm Span}(\psi_1, \ldots, \psi_I)$ is a good
approximation of ${\rm Ran}
\left(1_{[0,\beta^{-1/2}]}(-L^{(1)})\right)$. The support of $\psi_i$ is
essentially in a neighborhood of $z_i$ and Agmon estimates are used
  to prove exponential decay away from $z_i$.

The third step consists in projecting the approximation of $\nabla
v_1$ on the approximation of the eigenspace ${\rm Ran}
\left(1_{[0,\beta^{-1/2}]}(-L^{(1)})\right)$ using the following
result. Assume the following on the quasi-modes:
\begin{itemize}
\item {\em Normalization:} $\tilde v\in H^1_0({\rm e}^{-\beta V})$ and
  $\|\tilde v \|_{L^2({\rm e}^{-\beta V})} = 1$. For all $i \in \{1, \ldots,
  I\}$, $\psi_i \in  H^1_T({\rm e}^{-\beta V})$ and \footnote{The functional space
    $H^1_T({\rm e}^{-\beta V})$ is the space of $1$-forms in $H^1({\rm
      e}^{-\beta
      V})$  which satisfy the tangential Dirichlet boundary
    condition, see~\eqref{eq:L1}.}  $\|\psi_i\|_{L^2({\rm e}^{-\beta V})} = 1$.
 \item {\em Good quasi-modes:}
\begin{itemize}
\item $\forall \delta>0,  \|   \nabla  \tilde v\|_{L^2({\rm e}^{-\beta
      V})}^2       \     =  \    O ({\rm e}^{-\beta(V(z_1)-V(x_1) - \delta)}),$ 
\item      $\exists \varepsilon>0$, $\forall i \in \{1, \ldots, I\}$,
$  \|  1_{[\beta^{-1/2}, \infty)}(-L^{(1)})  \psi_i\|_{H^1({\rm e}^{-\beta
     V})}^2   \   =  \    O ({\rm e}^{-\beta (V(z_I) - V(z_1)  +     \varepsilon )}) $
 \end{itemize}
\item {\em Orthonormality of quasi-modes:}  $\exists \varepsilon_0>0$,
  $\forall i < j \in \{1, \ldots, I\}$,
$$\langle  \psi_i,  \psi_j\rangle_{L^2({\rm e}^{-\beta V})}=O (\ {\rm e}^{-\frac{\beta}{2}(V(z_j)- V(z_i)+\varepsilon_0)} \ ).$$
\item {\em Decomposition of $\nabla \tilde v$:} $\exists (C_i)_{1 \le
    i \le I} \in \R^I$, $\exists p>0$, $\forall i \in \{1, \ldots ,I\}$,
\begin{equation*}
  \langle       \nabla \tilde v  ,    \psi_i \rangle_{L^2({\rm e}^{-\beta V})}  =        C_i \ \beta^{-p}  {\rm e}^{-\frac{\beta}{2}(V(z_i)- V(x_1))}   \    (  1  +     O(\beta^{-1} ) \   ) .
\end{equation*}
\item {\em Normal components of the quasi-modes:} $\exists (B_i)_{1
    \le i \le I} \in \R^I$, $\exists m >0$, $\forall i,j \in \{1, \ldots ,I\}$,
\begin{equation*}
  \int_{\Sigma_i}     \psi_j \cdot n  \   {\rm e}^{- \beta V}    d\sigma =\begin{cases}   B_i \ \beta^{-m}   \     {\rm e}^{-\frac{\beta}{2} V(z_i)}  \    (  \   1 \  +     O(\beta^{-1} )  \   )   &  \text{ if } i=j,   \\
 0   &   \text{ if } i\neq j.
  \end{cases} 
  \end{equation*}
 \end{itemize}
 Then for $i=1,...,n$, when $\beta \to \infty$
\begin{equation*}
 \int_{\Sigma_i}   \   \partial_{n} v_1 \   {\rm e}^{-\beta V}
 d\sigma  =    C_i B_i  \ \beta^{-(p+m)}  \  {\rm e}^{-\frac{\beta}{2}(2V(z_i)- V(x_1))}    \   (1+  O(\beta^{-1}) ). 
  \end{equation*} 
The proof is based on a Gram-Schmidt orthonormalization
procedure. This result applied to the quasi-modes built in the second step yields~\eqref{eq:dnu1}.

\subsection{On the geometric assumption~\eqref{eq:agmon}}

In this section, we would like to discuss the geometric
assumption~\eqref{eq:agmon}. The question we would like to address is
the following: is such an assumption necessary to indeed prove the
result on the exit point density?

\begin{figure}[h]
\centering
  \includegraphics[height=4cm]{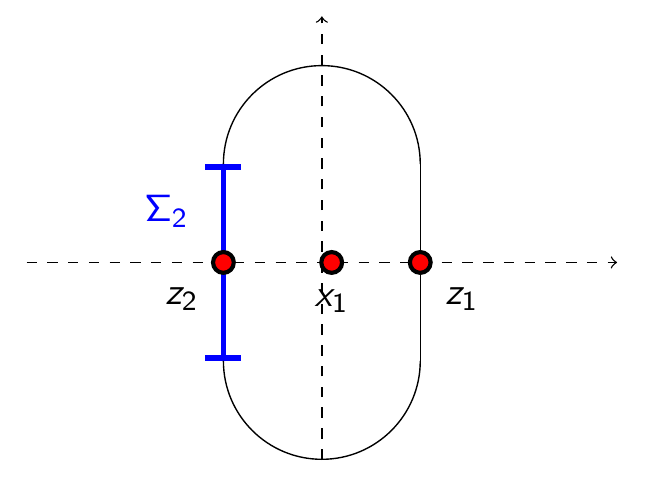}
  \caption{The domain $S$ is built as the union of the square with
    corners $(-1,-1)$ and $(1,1)$ and two half disks of radius $1$ and
  with centers $(0,1)$ and $(0,-1)$.}
  \label{fig:S_2d}
\end{figure}

In order to test this assumption numerically, we consider the
following simple two-dimensional setting. The potential function is
 $V(x,y)=x^2+y^2-ax$ with $a \in (0,1/9)$ on the domain $S$
 represented on Figure~\ref{fig:S_2d}. The two local minima on $\partial
 S$ are $z_1=(1,0)$ and $z_2=(-1,0)$. Notice that
 $V(z_2)-V(z_1)=2a>0$. The subset of the boundary around the highest
 saddle point is the segment $\Sigma_2$ joining the two points
 $(-1,-1)$ and $(-1,1)$. Using simple lower bounds on the Agmon distance,
 one can check that all the above assumptions are satisfied in this
 situation.

\begin{figure}[h]
\centering
\includegraphics[height=5cm]{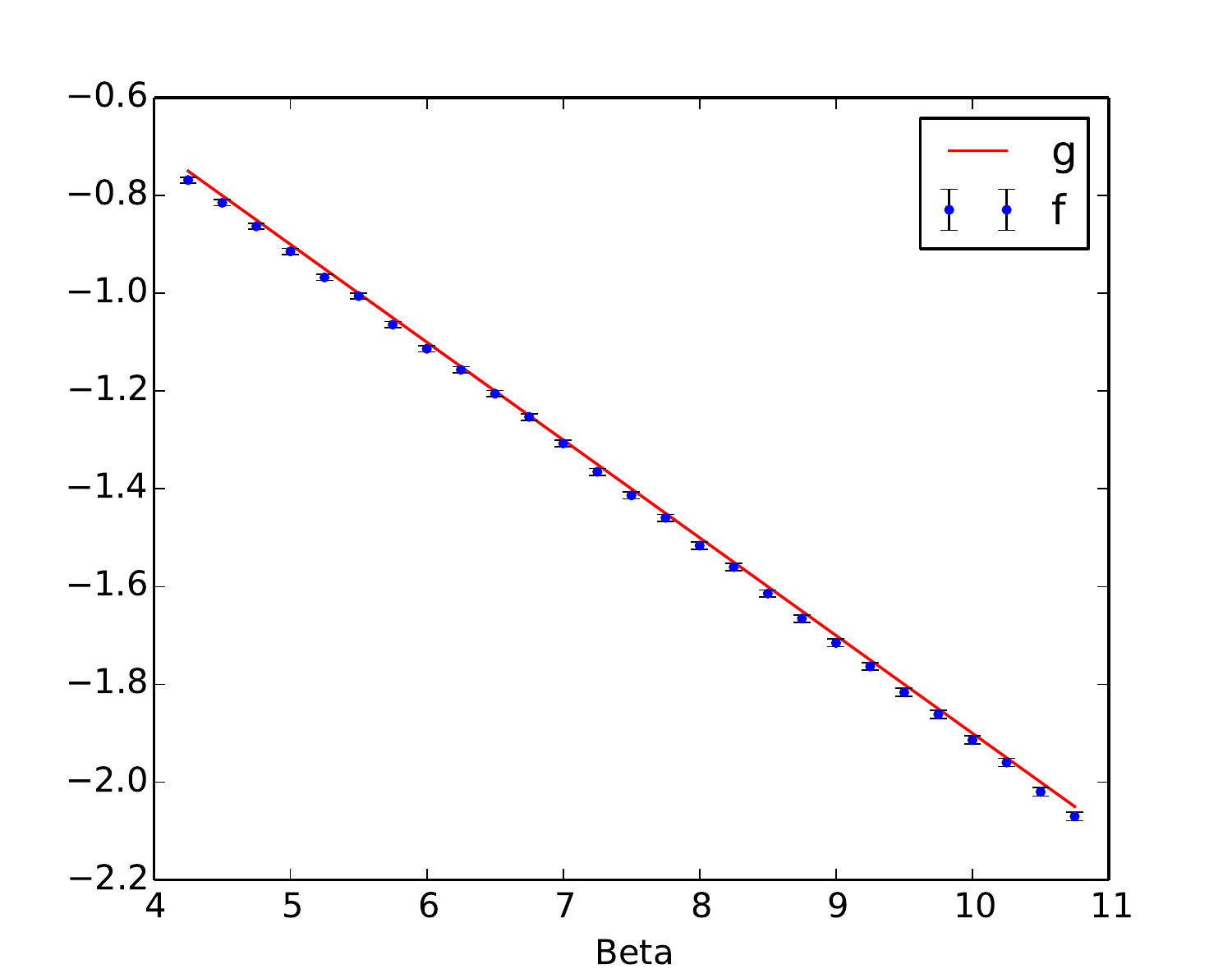}
  \caption{The probability $\P^{\nu_S}(X_{\tau_S} \in \Sigma_2)$:
 comparison of the theoretical result ($g$) with the numerical
    result ($f$, $\Delta t=5.10^{-3}$); $a=1/10$.}
  \label{fig:res1}
\end{figure}

\begin{figure}[h]
\centering
\includegraphics[height=5cm]{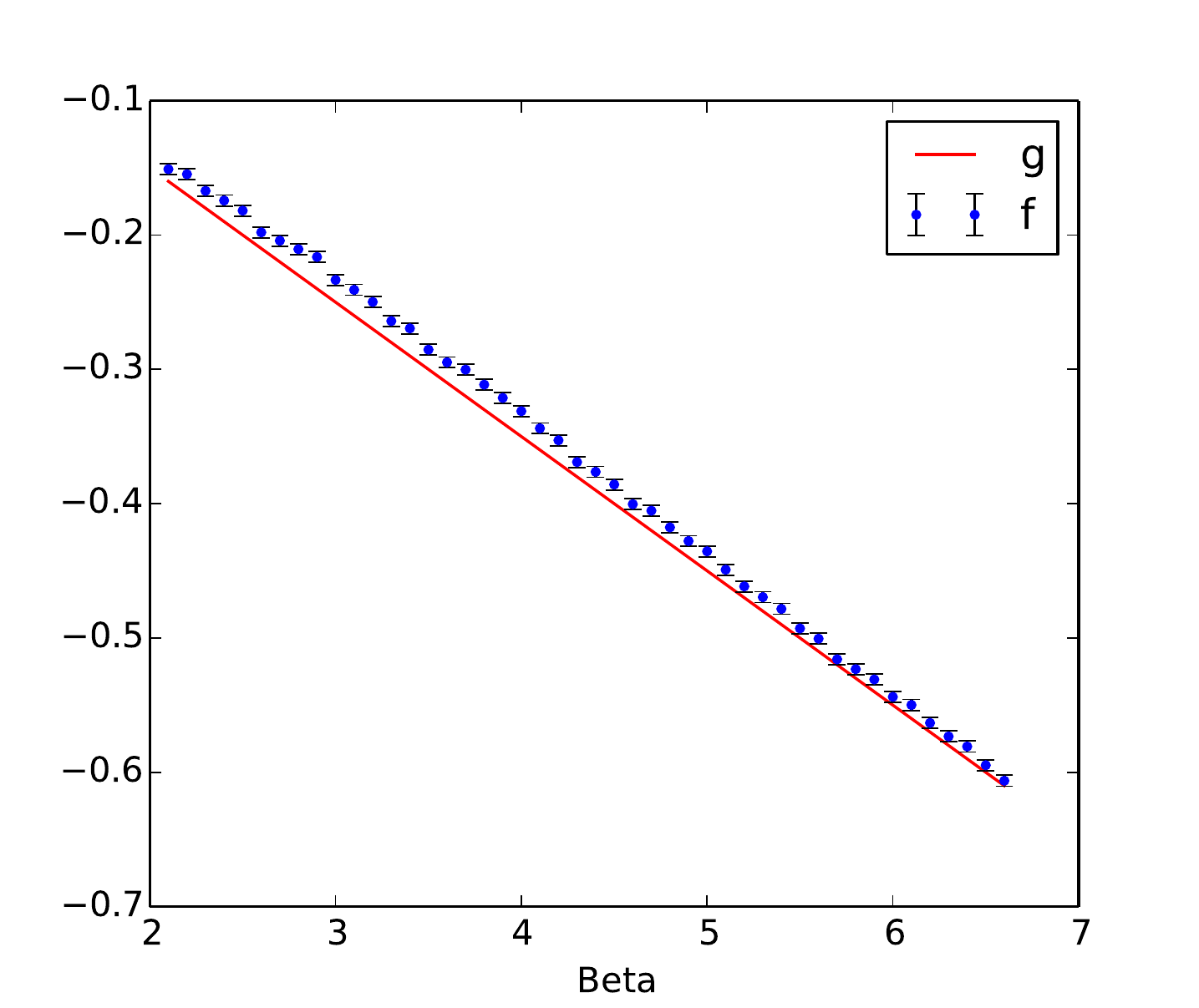}
  \caption{The probability $\P^{\nu_S}(X_{\tau_S} \in \Sigma_2)$:
    comparison of the theoretical result ($g$) with the numerical
    result ($f$, $\Delta t=2.10^{-3}$);  $a=1/20$.}
  \label{fig:res2}
\end{figure}

We then plot on Figures~\ref{fig:res1} ($a=1/10$) and~\ref{fig:res2}
($a=1/20$)  the numerically estimated probability
$f(\beta)=\P^{\nu_S}(X_{\tau_S} \in \Sigma_2)$, and compare it with
the theoretical result $g(\beta)=\frac{\partial_n V(z_2) \sqrt{ \det (\nabla^2 V|_{\partial S })
    (z_1) } }{\partial_n V(z_1) \sqrt{ \det (\nabla^2
    V|_{\partial S })   (z_2) }}   {\rm e}^{-\beta(V(z_2)-V(z_1))} $ (see
Equation~\eqref{eq:EK-p}). The probability $\P^{\nu_S}(X_{\tau_S} \in
\Sigma_2)$ is estimated using a Monte Carlo procedure, and the
dynamics~\eqref{eq:L} is discretized in time using an Euler-Maruyama
scheme with timestep $\Delta t$. We observe an excellent agreement
between the theory and the numerical results.

\begin{figure}[h]
\centering
\includegraphics[height=5cm]{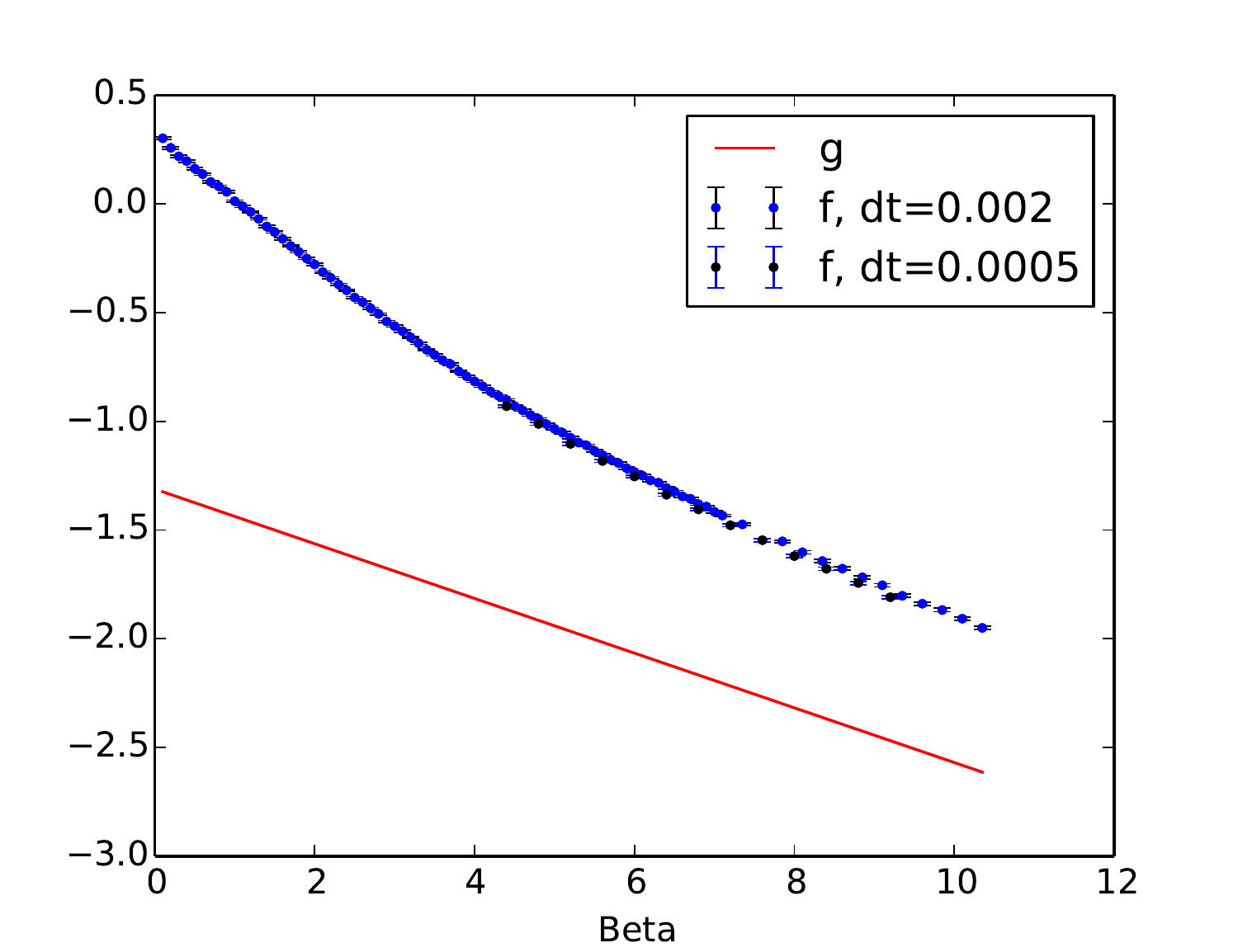}
  \caption{The probability $\P^{\nu_S}(X_{\tau_S} \in \Sigma_2)$:
      comparison of the theoretical result ($g$) with the numerical
    result ($f$, $\Delta t=2.10^{-3}$ and $\Delta t=5.10^{-4}$).}
  \label{fig:res3}
\end{figure}

Now, we modify the potential function $V$ in order not to satisfy
assumption~\eqref{eq:agmon} anymore. More precisely, the potential
function is $V(x,y)=(y^2-2 \ a(x) )^3$ with  $a(x)=a_1x^2+b_1x+ 0.5$
where $a_1$ and $b_1$ are chosen such that 
$a(-1+\delta)=0$, $a(1)=1/4$  for $\delta=0.05$. We have $V(z_1)=-1/8$ and
$V(z_2)=-8 (a(-1))^3> 0 > V(z_1)$. Moreover, two 'corniches' (which are in the level set $V^{-1}(\{0\})$ of $V$, and on which $|\nabla V|=0$) on
the 'slopes of the hills' of the potential $V$ join the point $(-1+\delta,0)$ to
$B_{z_2}^c$ so that $\inf \limits_{z\in B_{z_2}^c}
d_a(z,z_2) < V(z_2)-V(z_1)$ (the assumption~\eqref{eq:agmon} is not
satisfied). In addition  $V|_{\partial S}$  is a Morse function. The
function $V$ is not a Morse function on $S$, but an arbitrarily  small perturbation (which we neglect here)  turns it into a Morse function. When
comparing the numerically estimated probability
$f(\beta)=\P^{\nu_S}(X_{\tau_S} \in \Sigma_2)$, with
the theoretical result $g(\beta)$, we observe a discrepancy on the
prefactors, see Figure~\ref{fig:res3}.

Therefore, it seems that the construction of a jump Markov process
using the Eyring-Kramers law to estimate the rates to the neighboring
states is correct under some geometric assumptions. These geometric
assumptions appear in the proof when estimating rigorously the accuracy of the WKB
expansions as approximations of the quasi-modes of $L^{(1)}$.

\subsection{Concluding remarks}

In this section, we reported about some recent results obtained in~\cite{di-gesu-le-peutrec-lelievre-nectoux-16}.
We have shown that, under some geometric assumptions, the exit distribution from a state (namely the
law of the next visited state) predicted by a jump Markov process
built using the Eyring-Kramers formula is correct in the small
temperature regime, if the process starts from the QSD in the
state. We recall that this is a sensible assumption if the state is
metastable, and Equation~\eqref{eq:erreur} gives a
quantification of the error associated with this assumption.
Moreover, we have obtained bounds on the error introduced by
using the Eyring-Kramers formula. 

The analysis shows the importance of considering (possibly
generalized) saddle points on the boundary to identify the support of
the exit point distribution. This follows from the precise estimates we
obtain, which include the prefactor in the estimate of the probability
to exit through a given subset of the boundary. 

Finally, we checked by numerical experiments the fact that some
geometric assumptions are indeed required in order for all these
results to hold. These assumptions appear in the mathematical analysis when
rigorously justifying WKB expansions.

As mentioned above, we intend to generalize the results to the case
when the local minima of $V$ on $\partial S$ are saddle points.




\bibliography{ma_biblio,biblio_MD} 
\bibliographystyle{rsc} 

\end{document}